\documentclass[10pt]{amsart}
\usepackage{amssymb}
\usepackage{amscd}

\theoremstyle{definition}
\newtheorem{thm}{Theorem}[section]
\newtheorem{lem}[thm]{Lemma}
\newtheorem{prp}[thm]{Proposition}

\newtheorem{rmk}[thm]{Remark}
\newtheorem{ntn}[thm]{Notation}
\newtheorem{exa}[thm]{Example}

\newcommand{\beq}{\begin{equation}}
\newcommand{\eeq}{\end{equation}}
\newcommand{\beqr}{\begin{eqnarray*}}
\newcommand{\eeqr}{\end{eqnarray*}}
\newcommand{\bal}{\begin{align*}}
\newcommand{\eal}{\end{align*}}
\newcommand{\bei}{\begin{itemize}}
\newcommand{\eei}{\end{itemize}}

\newcommand{\af}{\alpha}
\newcommand{\bt}{\beta}
\newcommand{\gm}{\gamma}

\newcommand{\ep}{\varepsilon}

\newcommand{\et}{\eta}

\newcommand{\io}{\iota}
\newcommand{\te}{\theta}
\newcommand{\ld}{\lambda}
\newcommand{\sm}{\sigma}
\newcommand{\kp}{\kappa}
\newcommand{\ph}{\varphi}
\newcommand{\ps}{\psi}
\newcommand{\rh}{\rho}

\newcommand{\ta}{\tau}

\newcommand{\Gm}{\Gamma}

\newcommand{\Q}{{\mathbb{Q}}}
\newcommand{\Z}{{\mathbb{Z}}}
\newcommand{\R}{{\mathbb{R}}}
\newcommand{\C}{{\mathbb{C}}}
\newcommand{\N}{{\mathbb{N}}}
\newcommand{\F}{{\mathbb{F}}}

\newcommand{\id}{{\mathrm{id}}}

\newcommand{\Ker}{{\mathrm{Ker}}}

\newcommand{\Aut}{{\mathrm{Aut}}}
\newcommand{\spn}{{\mathrm{span}}}

\newcommand{\dirlim}{\displaystyle \lim_{\longrightarrow}}

\newcommand{\CO}{{\mathcal{O}}}

\newcommand{\Zq}[1]{\Z_{#1}}
\newcommand{\Zt}{\Zq{2}}
\newcommand{\Zqh}[1]{{\widehat{\Z}_{#1}}}
\newcommand{\Zth}{{\Zqh{2}}}

\newcommand{\Cs}[3]{{C^* (\Zq{#1}, {#2}, {#3})}}
\newcommand{\Cst}[2]{{\Cs{2}{#1}{#2}}}
\newcommand{\Csd}[3]{{C^* (\Zqh{#1}, {#2}, {#3})}}
\newcommand{\Cstd}[2]{{\Csd{2}{#1}{#2}}}

\newcommand{\oot}{of order two}
\newcommand{\ota}{order two automorphism}

\newcommand{\andeqn}{\,\,\,\,\,\, {\mbox{and}} \,\,\,\,\,\,}
\newcommand{\ts}[1]{{\textstyle{#1}}}

\newcommand{\ca}{C*-algebra}

\newcommand{\pj}{projection}

\newcommand{\hm}{homomorphism}

\newcommand{\ifo}{if and only if}

\newcommand{\uct}{Universal Coefficient Theorem}

\renewcommand{\S}{\subset}

\newcommand{\I}{\infty}

\title[Symmetries of Kirchberg algebras]{Symmetries
      of Kirchberg algebras}

\author{David J.\  Benson, Alex Kumjian, and N.\  Christopher Phillips}

\date{7 Feb.\  2003}

\address{Department of Mathematics, University of Georgia,
       Athens GA 30602-7403, USA.}

\email[]{djb@byrd.math.uga.edu}

\address{Department of Mathematics, University of Nevada,
       Reno NV 89557-0045, USA.}

\email[]{alex@unr.edu}

\address{Department of Mathematics, University  of Oregon,
       Eugene OR 97403-1222, USA.}

\email[]{ncp@darkwing.uoregon.edu}

\subjclass{Primary 20C10, 46L55;
 Secondary 19K99, 19L47, 46L40, 46L80.}
\thanks{Research partially supported by NSF grants
DMS-9988110 (D.\  J.\  Benson) and DMS-0070776 (N.\  C.\  Phillips).}

\begin{document}

\begin{abstract}
Let $G_0$ and $G_1$ be countable abelian groups.
Let $\gm_i$ be an automorphism of $G_i$ \oot.
Then there exists a unital Kirchberg algebra $A$ satisfying the \uct\  %
and with $[1_A] = 0$ in $K_0 (A)$,
and an automorphism $\af \in \Aut (A)$ \oot,
such that $K_0 (A) \cong G_0$, such that $K_1 (A) \cong G_1$, and
such that $\af_* \colon K_i (A) \to K_i (A)$ is $\gm_i$.
As a consequence, we prove that every $\Zt$-graded countable module
over the representation ring $R (\Zt)$ of $\Zt$
is isomorphic to the equivariant K-theory $K^{\Zt} (A)$ for some
action of $\Zt$ on a unital Kirchberg algebra $A$.

Along the way,
we prove that every not necessarily finitely generated
$\Z [ \Zt]$-module which is free as a $\Z$-module
has a direct sum decomposition with only three kinds of summands,
namely $\Z [ \Zt]$ itself
and $\Z$ on which the nontrivial element of $\Zt$ acts
either trivially or by multiplication by $-1$.
\end{abstract}

\maketitle

\section{Introduction}\label{Sec:Intro}

\indent
Following Definition~4.3.1 of Part~1 of~\cite{RoS}, we use the
term Kirchberg algebra for a purely infinite simple separable
nuclear \ca.
In this paper, we prove that any \ota\  of the
K-theory of a
unital Kirchberg algebra $A$ satisfying the \uct, and with
$[1_A] = 0$ in $K_0 (A)$, lifts to an \ota\  of $A$.

Recall the classification theorem for unital Kirchberg algebras
satisfying the \uct\  (\cite{Kr}; Theorem~4.2.4 of~\cite{Ph}):
if $A$ and $B$ are such algebras,
and if $\te \colon K_* (A) \to K_* (B)$ is a graded isomorphism
such that $\te ( [1_A]) = [1_B]$, then there is an
isomorphism $\ph \colon A \to B$ such that $\ph_* = \te$.
The statement also holds with ``isomorphism'' replaced by ``\hm''
everywhere (\cite{Kr}; in~\cite{Ph} see Theorem~4.1.3 and the proofs
of Corollary~4.4.2 and Theorem~4.2.4).
Moreover, every pair of countable abelian groups occurs as the
K-theory of a unital Kirchberg algebra satisfying the \uct, and
the $K_0$-class of the identity can be arbitrary
(Section~4.4 in Part~1 of~\cite{RoS} or Theorem~5.2 of~\cite{ER}).
Thus, the classification functor is surjective on
isomorphism classes in this case.

George Elliott has asked to what extent this functor ``splits'',
in the same sense
in which a surjection between abelian groups might (or might not) split.
That is: How close can one come to
constructing a functor $F$, from the values of the invariant, to
Kirchberg algebras satisfying the \uct,
such that $K_* (F (G_*)) \cong G_*$, etc.?
One can't do this exactly.
(See Example~\ref{NoSplit} below.)
We may then ask for some weaker sort of splitting,
or we can try to eliminate the problem by reducing the number
of morphisms in the categories.
The question we consider here, of lifting finite order
automorphisms of the K-theory to automorphisms of the algebra
of the same order, is a special case of what one gets by restricting the
morphisms to be isomorphisms.

Our method is to find, for given countable abelian groups $G_i$ and
\ota s $\gm_i$, some separable nuclear \ca\  $A_0$
whose K-theory is $G_*$ and some \ota\  $\af_0$
of $A_0$ which lifts the $\gm_i$.
In fact, $A_0$ will be type~I\@.
Then we apply the construction of~\cite{Kj} to get a Kirchberg
algebra with the same properties.
To construct $A_0$, we regard $G_*$ as a module over the group
ring $\Z [ \Zt]$.
(Throughout this paper, $\Z_p$ denotes $\Z / p \Z$,
not the $p$-adic integers.)
Generalizing an old result for the finitely generated case
(see Lemma~1 of~\cite{HR}),
we prove that every $\Z [ \Zt]$-module which is free as a $\Z$-module
has a direct sum decomposition with only three kinds of summands,
namely $\Z [ \Zt]$ itself
and $\Z$ on which the nontrivial element of $\Zt$ acts
either trivially or by multiplication by $-1$.
(Butler and Kov\'{a}cs have also obtained similar
results, by slightly different methods~\cite{BK};
a more general result is to appear in~\cite{BCK}.)
We produce $A_0$ by combining this structure theorem
with part of Schochet's geometric realization technique from~\cite{Sc}.

We want to explicitly point out that the K-theory of a \ca\  $A$
with an action of $\Zt$, regarded as a $\Z [ \Zt]$-module,
is not the same as the $\Zt$-equivariant K-theory of $A$,
regarded as a module over the representation ring
$R (\Zt) \cong \Z [ \Zt]$.
It is, however, the same as the $\Zt$-equivariant K-theory of
$C^* (\Zt, A)$, using the dual action.
(Compare with Lemma~\ref{CrProd}.)
Thus, as a corollary of our construction we obtain a realization
theorem, Theorem~\ref{RealizeEqKTh},
for $R (\Zt)$-modules as the equivariant K-theory of
actions of $\Zt$ on Kirchberg algebras.
There is also a version for type~I \ca s.

The first author would like to thank the MSRI for its hospitality
while this work was in progress.
The second and third authors would like to thank the organizers
of the conference
``Aperiodic Order, Dynamical Systems, Operator Algebras,
and Topology'', during which valuable discussions were held.
The third author would like to acknowledge the hospitality
of the Mathematisches Forschungsinstitut Oberwolfach,
where Elliott's raised his question and where
Example~\ref{M4O3} was constructed,
and of the University of Georgia for its invitation for a visit
in the spring of 2002.

This paper is organized as follows.
In the rest of this section, we demonstrate the failure
of the classification functor to split.
In Section~\ref{Sec:GpRMod} we prove a
structure theorem for $\Z [\Zt]$-modules
which are free as $\Z$-modules.
In Section~\ref{Sec:GRes} we prove the weak form of realization,
using some type~I \ca\  with the right K-theory.
In Section~\ref{Sec:PI} we prove the main result,
and also the realization theorems for $\Zt$-equivariant K-theory.
Section~\ref{Sec:ExplEx} contains an explicit formula
for the \ota\  in the simplest nontrivial case of our theorem,
namely the nontrivial automorphism of $K_* (M_3 (\CO_4))$.

Here is the example which shows that the most obvious sort of
splitting of the classification functor is not possible.

\begin{exa}\label{NoSplit}
There is no functor $F$, from the category of
$\Zt$-graded abelian groups $G_*$
with distinguished element $g$ in degree zero
and graded \hm s preserving
the distinguished elements, to the category of
unital Kirchberg algebras satisfying the \uct\  of~\cite{RS}
and unital \hm s,
such that $K_* (F (G_*, g)) \cong G_*$ with $[1_{F (G_*, g)}] \mapsto g$
for all $G_*$ and $g$, and such that the diagram
\[
\begin{CD}
K_* (F (G_*, g)) @>{F (\te)_*}>> K_* (F (H_*, h))\\
@VVV   @VVV\\
G_*  @>{\te}>> H_*
\end{CD}
\]
commutes whenever $\te \colon G_* \to H_*$ is a graded \hm\  such that
$\te (g) = h$.

Suppose we had such a functor.
Take
\[
G_0 = 0, \,\,\,\,\,\, G_1 = 0, \,\,\,\,\,\, H_0 = \Z,
\andeqn H_1 = 0.
\]
Take $g = 0$ and $h = 0$.
Let $\te \colon G_* \to H_*$ and $\rh \colon H_* \to G_*$
be the unique \hm s between these groups.
Then $\rh \circ \te = \id_{G_*}$.
Therefore, using functoriality,
\[
F (\rh) \circ F (\te) = F (\rh \circ \te) = F (\id_{G_*})
   = \id_{F (G_*)}.
\]
In particular, $F (\rh)$ is surjective.
However, $F (G_*, g)$
(which is isomorphic to the Cuntz algebra $\CO_2$~\cite{Cu1}) and
$F (H_*, h)$ (which is isomorphic to a suitable corner of $\CO_{\I}$)
are nonisomorphic simple \ca s, so there are no surjective \hm s from
$F (H_*, h)$ to $F (G_*, g)$.
\end{exa}

\begin{rmk}\label{NoSplitNonU}
The problem here has nothing to do with units, since the same
argument works even if one considers stable Kirchberg algebras
satisfying the \uct.
Moreover, one can even restrict, in either case, to the case
$K_1 = 0$.
\end{rmk}

\section{Modules over the group ring of $\Zt$}\label{Sec:GpRMod}

\indent
In this section, we prove a structure theorem for
modules over the group ring of $\Zt$ which are free as $\Z$-modules.
We have learned that Butler and Kov\'{a}c proved the main
result of this section a little earlier than we did,
using slightly different methods~\cite{BK}.
They use a sophisticated argument
to reduce to the countably generated case first,
which is avoided in our analysis.
There is also a subsequent generalization
to modules over the group ring of $\Zq{p}$
for an arbitrary prime $p$~\cite{BCK}.
We have decided to retain our proof for completeness of exposition,
and because it gives the result we need with a minimum of machinery.

For any unital ring $R$ and any discrete group $\Gm$, we let $R [\Gm]$
denote the ordinary algebraic group ring of $\Gm$ with coefficients
in $R$.
Although we are ultimately interested in the case $R = \Z$,
we will need the cases $R = \F_2$, the field with two elements,
and $R = \Q$.
Also, any $\Z [ \Gm ]$-module is automatically a $\Z$-module
(abelian group) by restriction of scalars.

The structure theorem is already known for the finitely generated case
(see the first lemma),
and the main point of this section is to remove
the finite generation hypothesis.

\begin{lem}\label{FGGpRMods}
Let $N$ be a finitely generated $\Z [ \Zt ]$-module which is
free as a $\Z$-module.
Then $N$ is a finite direct sum $\bigoplus_{i \in I} N_i$
in which each $N_i$ is isomorphic to one of the following three
$\Z [ \Zt ]$-modules:
\begin{itemize}
\item
$T_1 = \Z$ with the nontrivial element of $\Zt$ acting trivially.
\item
$T_2 = \Z$ with the nontrivial element of $\Zt$ acting
by multiplication by $-1$.
\item
$T_3 = \Z [ \Zt ]$.
\end{itemize}
The multiplicities of $T_1$, $T_2$ and $T_3$ in such a direct 
sum decomposition are independent of the decomposition.
\end{lem}

\begin{proof}
This is an immediate consequence of the canonical form for
invertible elements $a \in M_n (\Z)$ with $a^2 = 1$,
given in Lemma~1 of~\cite{HR},
and the discussion following the proof of that lemma.
For the matrix $L$ of~\cite{HR}, the computation
\[
L = \left( \begin{array}{cc} 1 & 0 \\ 1 & -1 \end{array} \right)
  = \left( \begin{array}{cc} 1 & -1 \\ 0 & 1 \end{array} \right)^{-1}
    \left( \begin{array}{cc} 0 & 1 \\ 1 & 0 \end{array} \right)
    \left( \begin{array}{cc} 1 & -1 \\ 0 & 1 \end{array} \right)
\]
shows that it is similar over $\Z$ to the
action of the nontrivial element of $\Zt$ on $\Z [ \Zt ]$ in its
usual basis.
\end{proof}

Theorem~74.3 of~\cite{CR} gives a structure theorem for
finitely generated $\Z [ \Zq{p} ]$-modules which are
free as a $\Z$-modules, for arbitrary primes $p$.
As described there,
the direct sum decomposition is in general not unique.

\begin{lem}\label{Freeness}
Let $M$ be a $\Z [ \Zt ]$-module.
Suppose $M$ is free as a $\Z$-module and
$M / 2 M$ is free as an $\F_2 [ \Zt ]$-module.
Then $M$ is free as a $\Z [ \Zt ]$-module.
\end{lem}

\begin{proof}
Theorem~6.1 of~\cite{BG}, with $G = \Zt$ and $H = \{ 0 \}$,
shows that $M$ is projective as a $\Z [ \Zt ]$-module.

If $M$ is finitely generated, apply Lemma~\ref{FGGpRMods}.
Since $T_1$ and $T_2$ are not projective as $\Z [ \Zt ]$-modules,
they can't appear in the direct sum.
Therefore $M$ is free.

So suppose $M$ is not finitely generated.
A direct calculation shows that $\Z [ \Zt ]$ has no
nontrivial idempotents.
(This is true for any finite group in place of $\Zt$,
by Corollary~8.1 of~\cite{Sw}.)
Therefore Corollary~4.5 of~\cite{Bs} applies,
and shows that $M$ is free.
\end{proof}

\begin{ntn}\label{Idem}
For the rest of this section, we
regard $\Z [\Zt]$ as a subring of $\Q [\Zt]$.
If $M$ is a $\Z [ \Zt ]$-module, then the corresponding
$\Q [ \Zt ]$-module is $\Q \otimes_{\Z} M$, and if $M$
is free as a $\Z$-module then we may clearly regard
$M$ as a submodule of $\Q \otimes_{\Z} M$.
We let $s$ be the nontrivial element of $\Zt$, and use
the same notation for the corresponding element of $\Z [\Zt]$.
Further set
\[
e = \ts{ \frac{1}{2} } ( 1 + s) \in \Q [\Zt] \andeqn
f = \ts{ \frac{1}{2} } ( 1 - s) \in \Q [\Zt],
\]
which are idempotents with $e + f = 1$.
\end{ntn}

\begin{lem}\label{L2}
Let $M$ be a $\Z [ \Zt ]$-module which is free as a $\Z$-module.
Let $m \in M$, and assume $(1 + s) m \in 2 M$.
Then, following Notation~\ref{Idem},
\[
m \in (M \cap e M) + (M \cap f M) \S \Q \otimes_{\Z} M.
\]
\end{lem}

\begin{proof}
We have $(1 - s) m = (1 + s) m - 2 m \in 2 M$.
Write $(1 + s) m = 2 m_0$ and $(1 - s) m = 2 m_1$
with $m_0, \, m_1 \in M$.
In $\Q \otimes_{\Z} M$ we have $m_0 = \ts{ \frac{1}{2} } ( 1 + s) = e m$
and $m_1 = f m$.
Then $m_0 \in M \cap e M$ and $m_1 \in M \cap f M$.
Moreover, $2 m = 2 (m_0 + m_1)$, so $m = m_0 + m_1$.
\end{proof}

\begin{lem}\label{L3}
Let $M$ be a $\Z [ \Zt ]$-module which is free as a $\Z$-module.
Following Notation~\ref{Idem}, suppose that
$M \cap e M = 2 e M$ and $M \cap f M = 2 f M$.
Then $M$ is free as a $\Z [ \Zt ]$-module.
\end{lem}

\begin{proof}
By Lemma~\ref{Freeness}, it suffices to show that
$M / 2 M$ is free as an $\F_2 [ \Zt ]$-module.

Set $N = (2 e M + 2 f M) / M$.
Consider $1 + s$ as a multiplication map on $M / 2 M$.
We claim that
\[
\Ker (1 + s) = (1 + s) (M / 2 M) = N.
\]

We consider $\Ker (1 + s)$ first.
Suppose $m \in M$ and $m + 2 M \in \Ker (1 + s)$.
Then $(1 + s) m \in 2 M$.
Using Lemma~\ref{L2} at the first step and the hypothesis at the
second, we get
\[
m \in (M \cap e M) + (M \cap f M) = 2 e M + 2 f M,
\]
so $m + 2 M \in N$.
For the reverse, suppose $m \in 2 e M + 2 f M$.
Write
\[
m = (1 + s) m_0 + (1 - s) m_1
\]
with $m_0, \, m_1 \in M$.
Then
\[
(1 + s) m = (1 + s)^2 m_0 + (1 + s) (1 - s) m_1
          = 2 (1 + s) m_0 \in 2 M,
\]
so $m + 2 M \in \Ker (1 + s)$.

Now we consider the range of $1 + s$.
Since $(1 + s)^2 = 0$ in $\F_2 [ \Zt ]$, we get
$(1 + s) (M / 2 M) \S \Ker (1 + s) = N$.
For the reverse, let $m \in 2 e M + 2 f M$, and again
write $m = (1 + s) m_0 + (1 - s) m_1$ with $m_0, \, m_1 \in M$.
Then
\begin{align*}
m + 2 M = (1 + s) (m_0 + m_1) - 2 s m_1 + 2 M
      & = (1 + s) (m_0 + m_1) + 2 M                 \\
      & \in (1 + s) (M / 2 M).
\end{align*}
This completes the proof of the claim.

Now $N$ is a $\F_2$-vector subspace of $M / 2 M$, so there exists
a $\F_2$-vector subspace $V \S M / 2 M$ such that
$V \oplus N = M / 2 M$.
We claim that $V \cap s V = \{ 0 \}$ and $V + s V = M / 2 M$.
To prove the first, suppose $v, \, w \in V$ and $v = s w$.
Then $(1 + s) v = (1 + s) s w = (1 + s) w$, so that
$(1 + s) (v - w) = 0$.
Therefore $v - w \in \Ker (1 + s) \cap V = N \cap V$, whence $v = w$.
Now $(1 + s) v = v + w = 2 v = 0$,
so $v \in (1 + s) (M / 2 M) \cap V = N \cap V$, whence $v = 0$.
For the second, let $d \in M / 2 M$.
Using the first claim, we have $N = (1 + s) (M / 2 M) = (1 + s) (V)$.
So we can write $d = x + (1 + s) w$ with $w, \, x \in V$.
Then $d = v + s w$ with $v = x + s w$ and $w$ both in $V$.
This proves the claim.

Given the last claim, it is immediate that
$M / 2 M \cong \F_2 [ \Zt ] \otimes_{\F_2} V$, which is a free
$\F_2 [ \Zt ]$-module.
\end{proof}

We only need the next two lemmas for the prime $2$.
Since they hold for arbitrary primes, with the same proof,
we may as well give them in that generality.

\begin{lem}\label{L4}
Let $M$ be a free $\Z$-module, and let $p$ be a prime number.
Let $V$ be an $\F_p$-vector subspace of $M / p M$.
Then there exists a direct summand $L$ of $M$ such that
the image of $L$ under the map $M \to M / p M$ is exactly $V$.
\end{lem}

\begin{proof}
Let $\pi \colon M \to M / p M$ be the quotient map.
Choose an $\F_p$-vector subspace $W$ of $M / p M$ such that
$M/ p M = V \oplus W$.
Let $\sm \colon M / p M \to W$ be the projection onto $W$ obtained
from this direct sum decomposition.
Using an $\F_p$-basis for $W$, find a free $\Z$-module $P_0$ with an
isomorphism $P_0 / p P_0 \to W$.
Let $\kp_0 \colon P_0 \to W$ be the composite of this isomorphism with
the quotient map $P_0 \to P_0 / p P_0$.
Since $M$ is a projective $\Z$-module and $\kp_0$ is surjective, there
exists a $\Z$-module \hm\  $\ph_0 \colon M \to P_0$ making the following
diagram commute:
\[
\begin{CD}
M @>\ph_0>> P_0\\
@V{\pi}VV   @VV{\kp_0}V\\
M / p M  @>\sm>> W
\end{CD}
\]

The image $\ph_0 (M)$ is a submodule of a free $\Z$-module, hence free
(Theorem~14.5 of~\cite{Fh}).
Let $( b_i )_{i \in I_0}$ be a $\Z$-basis for $\ph_0 (M)$.
Since $\kp_0 \circ \ph_0 = \sm \circ \pi$ is surjective,
$\{ \kp_0 (b_i) \colon i \in I_0 \}$ spans $W$.
Choose a subset $I \S I_0$ such that
$( \kp_0 (b_i) )_{ i \in I }$ is an $\F_p$-basis for $W$.
Set $P = \spn_{\Z} ( \{ b_i \colon i \in I \} ) \subset P_0$,
which is a free $\Z$-module with basis $( b_i )_{i \in I}$.
For each $j \in I_0 \setminus I$, choose $\af_{i, j} \in \F_p$
for $i \in I$, all but finitely many of which are zero,
such that $\kp_0 (b_j) = \sum_{i \in I} \af_{i, j} \kp_0 (b_i)$.
For $j \in I_0 \setminus I$ and $i \in I$, choose $\bt_{i, j} \in \Z$
whose image in $\F_p$ is $\af_{i, j}$, and such that if
$\af_{i, j} = 0$ then $\bt_{i, j} = 0$.
For each $j \in I_0 \setminus I$, all but finitely many of the
$\bt_{i, j}$ are zero.
Therefore there is a well defined surjective $\Z$-module
\hm\  $\mu \colon \ph_0 (M) \to P$ such that $p (b_j) = b_j$ for
$j \in I$ and $\mu (b_j) = \sum_{i \in I} \bt_{i, j} b_i$
for $j \in I_0 \setminus I$.
By construction, we have $\kp_0 \circ \mu (b_i) = \kp_0 (b_i)$ for all
$i \in I$.
Therefore $\left( \kp_0 |_P \right) \circ \mu = \kp_0$.
Setting $\ph = \mu \circ \ph_0$ and $\kp = \kp_0 |_P$, we thus
have a commutative diagram:
\[
\begin{CD}
M @>\ph>> P\\
@V{\pi}VV   @VV{\kp}V\\
M / p M  @>\sm>> W
\end{CD}
\]
The $\Z$-module $P$ is still free, and in addition $\ph$ is surjective.
Furthermore, $\Ker (\kp) = p P$, since a linear combination
of $( b_i )_{i \in I}$ has image in $W$ equal to zero
\ifo\  all coefficients are even.

Set $L = \Ker (\ph)$.
Since $P$ is projective,
there is a $\Z$-module \hm\  $\et \colon P \to M$
such that $\ph \circ \et = \id_P$.
Then $M = L \oplus \et (P)$.
We claim that $\pi (L) = V$.
First, if $m \in L$ then $\sm ( \pi (m)) = \kp (\ph (m)) = 0$,
so $\pi (m) \in \Ker (\sm) = V$.
Now suppose $v \in V$.
Choose $m_0 \in M$ such that $\pi (m_0) = v$.
Then $\kp (\ph (m_0)) = 0$, whence $\ph (m_0) \in p P$.
Write $\ph (m_0) = p x$ with $x \in P$.
Then $m = m_0 - p \et (x) \in \Ker (\ph)$ by a diagram chase,
and $\pi (m) = \pi (m_0) = v$.
This proves the claim, and the lemma.
\end{proof}

\begin{lem}\label{L5}
Let $M$ be a free $\Z$-module and let $p$ be a prime number.
Let $N$ be a $\Z$-submodule such that $p M \S N \S M$.
Then there exist $\Z$-submodules $N_0, \, N_1 \S M$
such that $M = N_0 \oplus N_1$ and $N = p N_0 \oplus N_1$.
\end{lem}

\begin{proof}
Let $\pi \colon M \to M / p M$ be the quotient map.
Using Lemma~\ref{L4}, write $M = N_0 \oplus N_1$
for $\Z$-submodules $N_0, \, N_1 \S M$, with $\pi (N_1) = \pi (N)$.
The result will now follow if we prove that $N = p N_0 + N_1$.

We prove $p N_0 + N_1 \S N$.
Clearly $p N_0 \S N$.
So let $m \in N_1$.
Choose $n \in N$ such that $\pi (n) = \pi (m)$.
Then $\pi (m - n) = 0$, so $m - n \in p M \S N$.
Thus $m \in N$.
We have $N_1 \S N$.

Now we prove the reverse inclusion.
Let $n \in N$.
Write $n = r + m_1$ with $r \in N_0$ and $m_1 \in N_1$.
Then $\pi (r) = \pi (n) - \pi (m_1)$.
Because $N_1 \S N$, we get $\pi (r) \in \pi (N_0) \cap \pi (N)$.
By considering a $\Z$-basis for $M$
made up of bases for $N_0$ and $N_1$,
we see that $\pi (N_0) \cap \pi (N) = \{ 0 \}$.
Therefore $\pi (r) = 0$.
So we can write $r = p m_0$ with $m_0 \in M$.
Since $M$ is torsion free and $N_0$ is a summand, we get $m_0 \in N_0$.
Therefore $n = p m_0 + m_1 \in p N_0 + N_1$.
\end{proof}

The uniqueness statement in the following theorem
will not be used in the rest of the paper,
but we include it because it seems of independent interest.

\begin{thm}\label{GpRingMods}
Let $M$ be a $\Z [ \Zt ]$-module which is
free as a $\Z$-module.
Then $M$ is a direct sum $\bigoplus_{i \in I} M_i$
in which each $M_i$ is isomorphic to one of the three
modules $T_1$, $T_2$, or $T_3$ of Lemma~\ref{FGGpRMods}.
The multiplicities of $T_1$, $T_2$ and $T_3$ in such a direct 
sum decomposition are independent of the decomposition.
\end{thm}

\begin{proof}
We use Notation~\ref{Idem}.

We first prove uniqueness of the multiplicities.
Suppose the $\Z [ \Zt ]$-module $M$ is isomorphic to the direct sum of
$n_1$ copies of $T_1$, $n_2$ copies of $T_2$,
and $n_3$ copies of $T_3$.
Then $n_1$ is the $\F_2$-vector space dimension of
$(M \cap e M) / 2 e M \S e M / 2 e M$,
and similarly $n_2 = \dim_{\F_2} ((M \cap f M) / 2 f M)$.
Finally, one checks that
$n_3 = \dim_{\F_2} (e M / (M \cap e M) )$,
which is the same as $\dim_{\F_2} (f M / (M \cap f M) )$.

Now we prove existence.
We have
\[
2 e M \S M \cap e M \S e M \andeqn 2 f M \S M \cap f M \S f M.
\]
Apply Lemma~\ref{L5} twice, to get $\Z$-submodules
$R_0, \, R_1 \S e M$ and $S_0, \, S_1 \S f M$ such that
\[
e M = R_0 \oplus R_1 \andeqn  M \cap e M = 2 R_0 \oplus R_1
\]
and
\[
f M = S_0 \oplus S_1 \andeqn M \cap f M = 2 S_0 \oplus S_1.
\]

Since $e M \cap f M = \{ 0 \}$ and $e + f = 1$, we can write
\[
(M \cap e M) \oplus (M \cap f M) \S M \S e M \oplus f M.
\]
Substituting from above, we have the inclusions of
$\Z$-module direct sums
\[
2 R_0 \oplus R_1 \oplus 2 S_0 \oplus S_1
  \S M \S R_0 \oplus R_1 \oplus S_0 \oplus S_1.
\]
Since $R_1, \, S_1 \S M$, we have the
$\Z$-module direct sum decomposition
\[
M = R_1 \oplus S_1 \oplus N \,\,\,\,\,\, {\mbox{with}} \,\,\,\,\,\,
N = M \cap (R_0 \oplus S_0).
\]
Now $R_1 \S e M$ and $s m = m$ for all $m \in e M$.
So $R_1$ is actually a $\Z [ \Zt ]$-module.
Since it is a $\Z$-submodule of $M$, it is a free $\Z$-module
(Theorem~14.5 of~\cite{Fh}),
and therefore as a $\Z [ \Zt ]$-module it is a direct
sum of $\Z [ \Zt ]$-modules isomorphic to $T_1$.
Similarly, from $S_1 \S f M$ and $s m = - m$ for all $m \in f M$,
we see that $S_1$ is a $\Z [ \Zt ]$-module which is a direct
sum of $\Z [ \Zt ]$-modules isomorphic to $T_2$.
We will complete the proof by showing that $N$ is a
$\Z [ \Zt ]$-module, so that the direct sum decomposition
above is really a $\Z [ \Zt ]$-module direct sum decomposition,
and that it is free, so that it is a direct
sum of $\Z [ \Zt ]$-modules isomorphic to $T_3$.

To prove that $N$ is a $\Z [ \Zt ]$-module, we simply
observe that $R_0$ and $S_0$ are $\Z [ \Zt ]$-modules,
for the same reason that $R_1$ and $S_1$ are.
For freeness, we know as above
that $N$ is a free $\Z$-module because it is a $\Z$-submodule of $M$.
By Lemma~\ref{L3}, it now suffices to prove that
$N \cap e N = 2 e N$ and $N \cap f N = 2 f N$.
That $2 e N \S N \cap e N$ and $2 f N \S N \cap f N$ is clear,
so we prove the reverse inclusions.

Let $m \in N \cap e N$.
Use $m \in N$ to write $m = a + b$ with $a \in R_0$ and $b \in S_0$.
Use $m \in e N$
to write $m = e (x + y)$ with $x \in R_0$ and $y \in S_0$.
Note that $e x = x$ and $e y = 0$,
and combine the two expressions for $m$ to
get $x - a = b$.
Since $x - a \in e M$ and $b \in f M$, we get $b = 0$.
Therefore $m = a \in R_0$.
{}From $m \in M \cap e M = 2 R_0 \oplus R_1$
we get $m \in 2 R_0$.
Now $R_0 \subset e M$ implies $e R_0 = R_0$,
so $R_0 \subset N$ implies $R_0 \subset e N$.
Therefore $m \in 2 e N$, as desired.

Now suppose $m \in N \cap f N$.
Write $m = a + b = f (x + y)$
with $a, \, x \in R_0$ and $b, \, y \in S_0$.
{}From $f x = 0$ and $f y = y$,
we get $y - b = a \in e M \cap f M = \{ 0 \}$.
So $m = b \in S_0$.
{}From $m \in M \cap f M = 2 S_0 \oplus S_1$
we now get $m \in 2 S_0$.
Since $S_0 \subset N \cap f M \subset f N$,
it follows that $m \in 2 f N$, as desired.
This completes the proof.
\end{proof}

\section{Geometric resolution}\label{Sec:GRes}

In this section, we prove that for every \ota\  $\gm$
of a countable $\Zt$-graded abelian group $G_*$,
there is some separable \ca\  $A$ with an \ota\  $\af \in \Aut (A)$,
such that there is a graded isomorphism $\mu \colon G_* \to K_* (A)$
which identifies $\gm$ with $\af_*$.
We borrow the method of geometric realization of resolutions
from Schochet's proof of the K\"{u}nneth formula
for \ca s~\cite{Sc}.
As it turns out, the \ca\  $A$ we produce will be of type~I
and nonunital.

Although we do not state this version formally,
the same proofs show that if $G_*$ is not countable,
we can still produce an automorphism of a nonseparable
type~I \ca\  satisfying the remaining conditions.

For any \ca\  $A$, we let $A^+$ denote its usual unitization,
in which we add a new unit even if $A$ is already unital.

\indent
\begin{lem}\label{ExistHom}
Let $C = \bigoplus_{n \in \N} [ C_0 ( \R) \oplus C_0 ( \R)]$,
and let $\gm \in \Aut (C)$ be the \ota\  given by
\[
\gm ( (f_1, g_1), \, (f_2, g_2), \, \dots)
  = ( (g_1, f_1), \, (g_2, f_2), \, \dots).
\]
Let $M$ be a subgroup of $K_1 (C)$ which is invariant under $\gm_*$.
Then there exists a separable Hilbert space $H$ with an action of
$\Zt$, a separable
commutative \ca\  $B$ with an automorphism $\bt$ \oot,
and a $\Zt$-equivariant \hm\  $\ph \colon B \to K (H) \otimes C$,
where $\Zt$ acts on $K (H)$ by conjugation,
such that
$K_0 (B) = 0$, such that $\ph_* \colon K_1 (B) \to K_1 (C)$ is
injective, and such that the image of $\ph_*$ is exactly $M$.
\end{lem}

\begin{proof}
The K-theory $K_* (C)$ is a countable abelian group with an
automorphism $\gm_*$ \oot, and hence in an obvious way
a $\Z [ \Zt ]$-module.
As $\Z [ \Zt ]$-modules,
$K_1 (C)$ is a countable direct sum of copies of $\Z [ \Zt ]$
and $K_0 (C) = 0$.
Moreover, the hypothesis on $M$ is just that it is a
$\Z [ \Zt ]$-submodule.
It is a free $\Z$-module because subgroups of free abelian groups
are free (Theorem~14.5 of~\cite{Fh}),
and it is clearly countably generated, so by
Theorem~\ref{GpRingMods} there is an isomorphism
$M \cong \bigoplus_{i \in I} M_i$ with $I$ finite or countable,
and with each $M_i$ isomorphic to one of the three modules $T_j$ of
Lemma~\ref{FGGpRMods}.
We regard each $M_i$ as a $\Z [ \Zt ]$-submodule of $K_1 (C)$.

For each $i \in I$, we construct a separable Hilbert space $H_i$
with an action of $\Zt$, a separable
commutative \ca\  $B_i$ with an automorphism $\bt_i$ \oot,
and a $\Zt$-equivariant
\hm\  $\ph_i \colon B_i \to K (H_i) \otimes C$, such that
$K_0 (B_i) = 0$, such that $(\ph_i)_* \colon K_1 (B_i) \to K_1 (C)$ is
injective, and such that the image of $(\ph_i)_*$ is exactly $M_i$.
There are three cases, but we begin by introducing notation common to
all three.
Define $u \colon \R \to S^1$ by
\[
u (t) = {\textstyle{ \exp \left( \pi i
    \left( 1 + t (1 + t^2)^{-1/2} \right) \right). }}
\]
Then $u$ is a unitary in the unitization $C_0 (\R)^+$ such that
$[u]$ generates $K_1 ( C_0 (\R))$, and $u - 1$ generates
$C_0 (\R)$ as a \ca.
In particular, if $E$ is any \ca\  and $v \in E^+$ is any unitary
such that $v - 1 \in E$, then there is a unique
\hm\  $\ps \colon C_0 (\R) \to E$ such that $\ps (u - 1) = v - 1$.
Further let $s$ be the nontrivial element of $\Zt$, and use the
same notation for the corresponding element of $\Z [ \Zt ]$.
We will also write elements of $K_1 (C)$ as sequences
\[
(k_1 + l_1 s, \, k_2 + l_2 s, \, \dots)
   \in \bigoplus_{n \in \N} \Z [ \Zt ],
\]
where, with all nontrivial entries being in the $n$-th position,
$(0, \dots, 0, 1, 0, \dots)$ corresponds to the $K_1$-class of
the unitary
\[
1 + ( (0, 0), \, \dots, \, (0, 0), \, (u - 1, \, 0),
        \, (0, 0), \, \dots) \in C^+
\]
and $(0, \dots, 0, s, 0, \dots)$ corresponds to the $K_1$-class of
the unitary
\[
1 + ( (0, 0), \, \dots, \, (0, 0), \, (0, \, u - 1),
        \, (0, 0), \, \dots) \in C^+.
\]

The first case is $M_i \cong T_1$.
Let $m$ be a generator of $M_i$ as a $\Z$-module.
Then $s m = m$.
Therefore $m$ has the form
\[
m = (k_1 + k_1 s, \, k_2 + k_2 s, \, \dots),
\]
with all but finitely many of the $k_j$ equal to zero.
Take $H_i = \C$ with the trivial action of $\Zt$,
take $B_i = C_0 (\R)$, take $\bt_i = \id_{B_i}$, and take
$\ph_i$ to be the \hm\  determined by
\[
\ph_i (u - 1) =
 ( (u^{k_1} - 1, \, u^{k_1} - 1), \, (u^{k_2} - 1, \, u^{k_2} - 1), \,
   \dots).
\]
The required properties are immediate.

Next, suppose $M_i \cong T_2$.
Let $m$ be a generator of $M_i$ as a $\Z$-module.
Then $s m = - m$.
Therefore $m$ has the form
\[
m = (k_1 - k_1 s, \, k_2 - k_2 s, \, \dots),
\]
with all but finitely many of the $k_j$ equal to zero.
Take $H_i = \C$ with the trivial action of $\Zt$, and
take $B_i = C_0 (\R)$.
Take $\bt_i \colon C_0 (\R) \to C_0 (\R)$
to be the \hm\  $\bt_i (f)(t) = f (- t)$,
which is the unique \hm\  such that $\bt_i (u - 1) = u^* - 1$.
Take $\ph_i$ to be the \hm\  determined by
\[
\ph_i (u - 1) =
 ( (u^{k_1} - 1, \, u^{- k_1} - 1), \, (u^{k_2} - 1, \, u^{- k_2} - 1),
    \, \dots).
\]
Equivariance follows from
\[
\ph_i (u^* - 1) =
 ( (u^{- k_1} - 1, \, u^{k_1} - 1), \, (u^{- k_2} - 1, \, u^{k_2} - 1),
    \, \dots),
\]
and the rest of the required properties are immediate.

Finally, suppose $M_i \cong T_3$.
Let $m \in K_1 (C)$ be the image of $1 \in T_3$ under the
isomorphism $T_3 \to M_i$.
Write $m = (k_1 + l_1 s, \, k_2 + l_2 s, \, \dots)$.
Then the image of $s$ is
\[
s m = (l_1 + k_1 s, \, l_2 + k_2 s, \, \dots),
\]
and these two elements are linearly independent over $\Z$.
Take $H_i = \C^2$, with $s ( \xi_1, \xi_2) = ( \xi_2, \xi_1)$.
Take $B_i = C_0 (\R) \oplus C_0 (\R)$,
with $\bt_i (f_1, f_2) = (f_2, f_1)$.
To define
\[
\ph_i \colon C_0 (\R) \oplus C_0 (\R) \to L (\C^2) \otimes C,
\]
let $p_1, \, p_2 \in L (\C^2)$ be the projections on the
first and second coordinates, and let
$\ps_1, \, \ps_2 \colon C_0 (\R) \to C$ be the unique \hm s such
that
\[
\ps_1 (u - 1) =
 ( (u^{k_1} - 1, \, u^{l_1} - 1), \, (u^{k_2} - 1, \, u^{l_2} - 1),
    \, \dots)
\]
and
\[
\ps_2 (u - 1) =
 ( (u^{l_1} - 1, \, u^{k_1} - 1), \, (u^{l_2} - 1, \, u^{k_2} - 1),
    \, \dots).
\]
Then define
\[
\ph_i (f_1, f_2) = p_1 \otimes \ps_1 (f_1) + p_2 \otimes \ps_2 (f_2).
\]
This is a \hm\  because $p_1$ and $p_2$ are orthogonal.
Equivariance follows from the fact that the action of $s$
exchanges $p_1$ and $p_2$, and the formulas for $m$ and $s m$.
That $(\ps_i)_*$ is injective follows from the fact that
$m$ and $s m$ are linearly independent over $\Z$,
and the remaining required properties are clear.

Now let $H$ be the Hilbert direct sum $H = \bigoplus_{i \in I} H_i$,
and let $B$ be the \ca\  direct sum $B = \bigoplus_{i \in I} B_i$.
Give both the action of $\Zt$ coming from the actions on the
summands.
Define $\ph \colon B \to K (H) \otimes C$ by taking
$\ph ( (b_i)_{i \in I} )$ to be the block diagonal element
$\bigoplus_{i \in I} \ph_i (b_i)$.
This element is in fact in $K (H) \otimes C$ because for every
$\ep > 0$, we have $\| b_i \| < \ep$ for all but finitely many $i$.
Then $\ph$ is the required \hm.
\end{proof}

\begin{lem}\label{ExistAlg0}
Let $G$ be a countable abelian group, and let $\nu \colon G \to G$
be an automorphism of $G$ \oot.
Then there exist a separable type~I \ca\  $A$
such that $K_1 (A) = 0$,
an automorphism $\af \in \Aut (A)$ \oot,
and an isomorphism $\mu \colon G \to K_0 (A)$,
such that $\mu \circ \nu = \af_* \circ \mu$.
\end{lem}

\begin{proof}
The group $G$ is a $\Z [ \Zt ]$-module in an obvious way,
and it is countable by assumption.
Set $N = \bigoplus_{n \in \N} \Z [ \Zt ]$, and
choose a surjective $\Z [ \Zt ]$-module \hm\  %
$\ta \colon N \to G$.
Set $M = \Ker (\ta)$.
Let $C$ be as in Lemma~\ref{ExistHom},
and identify $N$ with $K_1 (C)$.
Apply Lemma~\ref{ExistHom} with this $M$, obtaining
a $\Zt$-equivariant \hm\  of separable type~I \ca s
$\ph \colon B \to K (H) \otimes C$
such that (by abuse of notation)
$K_1 (B) = M$, $K_1 ( K (H) \otimes C ) = N$, and
$\ph_*$ is the inclusion.
To simplify the notation, set $D = K (H) \otimes C$.

Let $A$ be the mapping cone
\[
A = \{ (f, b) \in C ([0, 1], \, D) \oplus B \colon
  {\mbox{$f (0) = 0$ and $f (1) = \ph (b)$}} \}.
\]
Since $\ph$ is equivariant, this algebra has an obvious action of
$\Zt$.
(The action is trivial on $[0, 1]$.)
Moreover, with $S D = C_0 ( (0, 1), \, D)$ being the usual
suspension of $D$, there is an equivariant short exact sequence
\[
0 \longrightarrow S D \longrightarrow
A \longrightarrow B \longrightarrow 0.
\]
By naturality the corresponding six term exact sequence in
K-theory is actually a sequence of $\Z [ \Zt ]$-modules.
Since $K_0 (B) = 0$ and $K_0 (D) = 0$, this sequence reduces to
\[
0 \longrightarrow K_1 (A) \longrightarrow
K_1 (B) \longrightarrow K_0 (S D) \longrightarrow
K_0 (A) \longrightarrow 0.
\]
Using the Bott periodicity isomorphism $K_0 (SD) \cong K_1 (D)$,
and appropriately identifying the maps following
Theorem~3.5 of~\cite{Sc3},
this sequence can be naturally identified with 
\[
0 \longrightarrow K_1 (A) \longrightarrow
K_1 (B) \stackrel{\ph_*}{\longrightarrow} K_1 (D) \longrightarrow
K_0 (A) \longrightarrow 0.
\]
By our construction and by naturality of the sequence,
we therefore have
the sequence of $\Z [ \Zt ]$-modules
\[
0 \longrightarrow K_1 (A) \longrightarrow
M \stackrel{\ta}{\longrightarrow} N \longrightarrow
K_0 (A) \longrightarrow 0,
\]
from which it follows that $K_1 (A) = 0$ and that
$K_0 (A) \cong N / M \cong G$ as $\Z [ \Zt ]$-modules.
\end{proof}

\begin{prp}\label{ExistAlg}
Let $G_0$ and $G_1$ be countable abelian groups.
Let $\gm_i$ be an automorphism of $G_i$ \oot.
Then there exist a separable type~I \ca\  $A$,
an automorphism $\af \in \Aut (A)$ \oot,
and a graded isomorphism $\mu \colon G_* \to K_* (A)$,
such that the diagram
\[
\begin{CD}
G_*  @>{\gm_*}>> G_*\\
@V{\mu}VV   @VV{\mu}V\\
K_* (A) @>{\af_*}>> K_* (A)
\end{CD}
\]
commutes.
\end{prp}

\begin{proof}
Apply Lemma~\ref{ExistAlg0} to $G_0$ and $\gm_0$,
obtaining a \ca\  $A_0$ with an \ota\  $\af_0$.
Apply Lemma~\ref{ExistAlg0} to $G_1$ and $\gm_1$,
obtaining a \ca\  $A_1$ with an \ota\  $\af_1$.
Take $A = A_0 \oplus S A_1$ and $\af = \af_0 \oplus S \af_1$.
This \ca\ and automorphism satisfy the conclusions by Bott periodicity.
\end{proof}

\section{From type~I to purely infinite simple}\label{Sec:PI}

\indent
For countable abelian groups $G_0$ and $G_1$
and \ota s $\gm_i$ of $G_i$,
we can now produce a separable type~I C*-algebra $A$ with K-theory
$G_0 \oplus G_1$, and an \ota\  which induces
$\gm_i$ on K-theory.
We want to use the construction of~\cite{Kj}
to produce a unital Kirchberg algebra.
However, $A$ is not unital, and this construction
will then not produce a unital \ca.
Unitizing $A$ changes the K-theory.
To remedy this problem, we introduce in the next proposition a
kind of unitization functor which does not change the K-theory.
Our functor preserves nuclearity and the \uct,
but not type~I (or simplicity or finiteness).
Although we will not use this fact, it also seems interesting to
observe that the functor sends exact \ca s to exact \ca s .

\begin{prp}\label{UnitFunc}
There is a functor $F$ from the category of all \ca s and
\ca\  \hm s to the subcategory of all unital \ca s and
unital \ca\  \hm s, and a natural transformation $\et$ from
the identity functor to $F$, such that:
\begin{itemize}
\item[(1)]
For every \ca\  $A$, the \hm\  $\et_A \colon A \to F (A)$ is
a KK-equivalence, and in particular is an isomorphism on K-theory.
\item[(2)]
If $A$ is separable, or nuclear, or exact, or satisfies the
\uct\  of~\cite{RS}, then
$F (A)$ also has the same property.
\end{itemize}
\end{prp}

\begin{proof}
Choose and fix a nonzero \pj\  $e \in \CO_{\I}$ such that
$[e] = 0$ in $K_0 (\CO_{\I})$, a \pj\  $q \leq e$ such that
$[q] = [1_{\CO_{\I} } ]$ in $K_0 (\CO_{\I})$, and a unital subalgebra
$D \S e \CO_{\I} e$ such that $D \cong \CO_2$.
For any \ca\  $A$, let $A^+$ be the unitization
(as in Section~\ref{Sec:GRes}),
let $\pi_A \colon A^+ \to \C$ be the canonical map, and define
\[
F (A) = \{ b \in e \CO_{\I} e \otimes A^+ \colon
  ( \id_{e \CO_{\I} e} \otimes \pi_A) (b) \in D \}.
\]
Let $\io_A \colon A \to A^+$ be the inclusion, and define
$\et_A \colon A \to F (A)$ by $\et_A (a) = q \otimes \io_A (a)$.
This element is in $F (A)$ because
$( \id_{e \CO_{\I} e} \otimes \pi_A) (q \otimes \io_A (a)) = 0$.

Since the map $\ld \mapsto \ld q$ from $\C$ to $e \CO_{\I} e$
is a KK-equivalence, so is the map $a \mapsto q \otimes a$
from $A$ to $e \CO_{\I} e$.
(See Example~19.1.2(c) of~\cite{Bl}.)
Moreover, there is a split short exact sequence
\[
0 \longrightarrow e \CO_{\I} e \otimes A
  \longrightarrow F (A) \longrightarrow \CO_2 \longrightarrow 0.
\]
Since $\CO_2$ is KK-equivalent to the zero \ca,
it follows that, for every separable \ca\  $B$, the inclusion
$\ph \colon e \CO_{\I} e \otimes A \to F (A)$ induces isomorphisms
\[
\ph_B^* \colon KK (F (A), \, B) \to KK (e \CO_{\I} e \otimes A, \, B)
\]
and
\[
\ph^B_* \colon KK (B, \, e \CO_{\I} e \otimes A) \to KK (B, \, F (A)).
\]
It now follows from the Yoneda Lemma that $\ph$ is a KK-equivalence
(see Section~III.2 of~\cite{Mc}),
but we give the easy direct argument here.
In terms of the Kasparov product,
these isomorphisms have the formulas $\ph_B^* (\gm) = [\ph] \cdot \gm$
and $\ph^B_* (\gm) = \gm \cdot [\ph]$.
(See Proposition~18.7.2 of~\cite{Bl}.)
Choose
\[
\af \in KK (F (A), \, e \CO_{\I} e \otimes A)
\andeqn
\bt \in KK (e \CO_{\I} e \otimes A, \, F (A))
\]
such that
\[
\ph_{e \CO_{\I} e \otimes A}^* (\af) = [ \id_{e \CO_{\I} e \otimes A}]
\andeqn
\ph^{F (A)}_* (\bt) = [\id_{F (A)}].
\]
Then
\[
[\ph] \cdot \af = [ \id_{e \CO_{\I} e \otimes A}]
\andeqn
\bt \cdot [\ph] = [\id_{F (A)}],
\]
from which it follows that $\ph$ is also a KK-equivalence.
So $\et_A$ is a KK-equivalence.
This proves the property~(1).

It is immediate from the KK-equivalence of~(1)
that $F (A)$ satisfies the \uct\  if $A$ does.
That $F$ preserves separability and nuclearity
is immediate from the split exact sequence above, since
$e \CO_{\I} e$ and $\CO_2$ are separable and nuclear.
The same argument works for exactness:
the minimal tensor product of exact \ca s is
exact by Proposition~7.1(iii) of~\cite{Kr0},
and extensions with completely positive
splittings preserve exactness by Proposition~7.1(vi) of~\cite{Kr0}.
\end{proof}


Next, we need an equivariant version of the construction of~\cite{Kj}.

\begin{prp}\label{MakePI}
Let $\af \colon \Gm \to \Aut (A)$ be an action of a countable discrete
group $\Gm$ on a separable unital \ca\  $A$.
Then there exists a separable unital purely infinite simple \ca\  $B$,
an action $\bt \colon \Gm \to \Aut (B)$, and an equivariant 
unital \hm\  $\ph \colon A \to B$, such that $\ph$ is a KK-equivalence.
Moreover, if $A$ is nuclear then so is $B$, and if $A$ satisfies
the \uct\  then so does $B$.
\end{prp}

\begin{proof}
Choose an injective unital representation $\pi_0 \colon A \to L (H_0)$
of $A$ on a separable Hilbert space $H_0$ such that
$\pi_0 (A) \cap K (H_0) = \{ 0 \}$.
Let $H = l^2 (\Gm, H_0)$, and let
$u \colon \Gm \to U (H)$ and $\pi \colon A \to L (H)$ be the
components of the regular covariant representation of $(\Gm, A)$
associated to $\pi_0$, as in~7.7.1 in~\cite{Pd}.
Then $u_{\gm} \pi (a) u_{\gm}^* = \pi ( \af_{\gm} (a))$ for all
$a \in A$ and $\gm \in \Gm$.
(See~\cite{Pd}.)
Moreover, $\pi (A) \cap K (H) = \{ 0 \}$.

Now we follow the proof of Proposition~2.1 of~\cite{Kj}.
Let $E = H \otimes_{\C} A$ be the Hilbert $A$-bimodule defined
there, and let $\ph \colon A \to L (E)$ be as there. 
Each $\gm \in \Gm$ induces a isometric $\C$-linear map
$v_{\gm} \colon E \to E$
given by $v_{\gm} (\xi \otimes a) = u_{\gm} \xi \otimes \af_{\gm} (a)$
for $\xi \in H$ and $a \in A$.
These maps have the following properties:
\begin{itemize}
\item[(1)]
$v_{\gm} ( \xi a) = v_{\gm} ( \xi ) \af_{\gm} (a)$
for $\gm \in \Gm$, $\xi \in E$, and $a \in A$.
\item[(2)]
$\langle v_{\gm} \xi, \, v_{\gm} \et \rangle
   = \af_{\gm} (\langle \xi, \et \rangle)$
for $\gm \in \Gm$ and $\xi, \, \et \in E$.
\item[(3)]
$\| v_{\gm} \xi \| = \| \xi \|$ for $\gm \in \Gm$ and $\xi \in E$.
\item[(4)]
$v_{\gm^{-1}} = v_{\gm}^{-1}$ for $\gm \in \Gm$.
\item[(5)]
$v_{\gm} \ph (a) v_{\gm}^{-1} = \ph ( \af_{\gm} (a))$
for $\gm \in \Gm$ and $a \in A$.
\end{itemize}
Caution: $v_{\gm}$ is not a right $A$-module \hm.
Nevertheless, these properties imply that
$b \mapsto v_{\gm} b v_{\gm}^{-1}$ is a *-automorphism
of the \ca\  $L (E)$ of all adjointable right $A$-module
morphisms of $E$;
moreover, $\ph$ is equivariant.

Next, let ${\mathcal{E}}_+ = \bigoplus_{n = 0}^{\infty} E^{\otimes n}$
be the Fock space of $E$, as after Definition~1.3 of~\cite{Kj}
or as at the beginning of Section~1 of~\cite{Pm}.
We let $\gm \in \Gm$ act on $E^{\otimes n}$ via $v_{\gm}^{\otimes n}$
(via $\af_{\gm}$ on $E^{\otimes 0} = A$),
and on ${\mathcal{E}}_+$ via the direct sum $w_{\gm}$ of these actions.
This action is still isometric, and we again have the
analogs of Properties (1)--(5) above.
The role of $\ph$ is now played by the
\hm\  $\ph_+ \colon A \to L ({\mathcal{E}}_+)$ of the
discussion after Definition~1.3 of~\cite{Kj}.
Following the construction there, the Toeplitz algebra
${\mathcal{T}}_E$ is by definition the C*-subalgebra of
$L ({\mathcal{E}}_+)$ generated by the operators $T_{\xi}$ defined
there, for $\xi \in E$.
One checks that $T_{v_{\gm} \xi} = w_{\gm} T_{\xi} w_{\gm}^{-1}$
for $\gm \in \Gm$ and $\xi \in E$.
Therefore $\Gm$ acts on ${\mathcal{T}}_E$ via *-automorphisms.
Since $\pi (A) \cap K (H) = \{ 0 \}$, Lemma~2.1 of~\cite{Kj} and
Corollary~3.14 and Theorem~3.13 of~\cite{Pm} show that the
canonical map from ${\mathcal{T}}_E$ to the
Cuntz-Pimsner algebra ${\mathcal{O}}_E$ is an isomorphism.
Therefore we have an action of $\Gm$ on ${\mathcal{O}}_E$,
and $\ph_+ \colon A \to {\mathcal{O}}_E$ is equivariant.

Theorem~2.8 of~\cite{Kj} shows that
${\mathcal{O}}_E$ is purely infinite and simple.
Separability is clear.
The proof of Theorem~3.1 of~\cite{Kj} shows that if $A$ is
nuclear, then so is ${\mathcal{O}}_E$.
Corollary~4.5 of~\cite{Pm} shows that $\ph_+$ is a KK-equivalence,
from which it is immediate that if $A$ satisfies the \uct\  then so
does ${\mathcal{O}}_E$.
\end{proof}

\begin{thm}\label{LiftGpAction}
Let $G_0$ and $G_1$ be countable abelian groups.
Let $\gm_i$ be an automorphism of $G_i$ \oot.
Then there exist a unital Kirchberg algebra $A$ satisfying the \uct\  %
and with $[1_A] = 0$ in $K_0 (A)$,
an automorphism $\af \in \Aut (A)$ \oot,
and a graded isomorphism $\mu \colon G_* \to K_* (A)$,
such that the diagram
\[
\begin{CD}
G_*  @>{\gm_*}>> G_*\\
@V{\mu}VV   @VV{\mu}V\\
K_* (A) @>{\af_*}>> K_* (A)
\end{CD}
\]
commutes.
\end{thm}

\begin{proof}
Use Proposition~\ref{ExistAlg} to find a separable type~I \ca\  $A_0$,
an automorphism $\af_0 \in \Aut (A_0)$ \oot,
and a graded isomorphism $\mu_0 \colon G_* \to K_* (A_0)$,
such that $\mu_0 \circ \gm = (\af_0)_* \circ \mu_0$.

Let $F$ be the functor of Proposition~\ref{UnitFunc},
and let $\et \colon A_0 \to F (A_0)$ be the natural transformation from
there.
By Part~(2) of Proposition~\ref{UnitFunc},
the algebra $F (A_0)$ is separable, unital, nuclear, and satisfies the \uct;
also, $\et$ is a KK-equivalence.
Since $F$ is a functor, $F (\af_0)$ is an \ota\  %
of $F (A_0)$, and naturality implies that
$\et \circ \af_0 = F (\af_0) \circ \et$.

Now use Proposition~\ref{MakePI} to find a
unital Kirchberg algebra $A$ satisfying the \uct,
an automorphism $\af \in \Aut (A)$ \oot,
and a $\Zt$-equivariant \hm\  $\ph \colon F (A_0) \to A$
which is a KK-equivalence.
Then $\ph \circ \et \colon A_0 \to A$ is $\Zt$-equivariant
and is a KK-equivalence.
The theorem is therefore proved by taking
$\mu = \ph_* \circ \et_* \circ \mu_0$.
\end{proof}

It should be easy to arrange to have
$[1_A]$ correspond to any element in $G$ of the form
$\et + \gm_0 (\et)$, but we don't know how to get arbitrary
$\gm_0$-invariant elements of $G_0$.

We now turn to the problem of realizing $R (\Zt)$-modules as
equivariant K-theory for actions on \ca s.
The main point is contained in the next lemma,
which we state in greater generality.

Before stating it, recall (Section~2.2 of~\cite{PhT})
that if a compact group $G$ acts on a \ca\  $A$,
then the equivariant K-theory $K_*^G (A)$ is, in a canonical way,
a module over the representation ring $R (G)$ of $G$.
Further recall that when $G$ is abelian,
the representation ring $R (G)$ is just $\Z [ {\widehat{G}} ]$;
in particular, if $G$ is a countable abelian group, then
$R ({\widehat{G}})$ is canonically isomorphic to $\Z [G]$.

\begin{lem}\label{CrProd}
Let $G$ be a countable abelian group,
and let $\af \colon G \to \Aut (A)$
be an action of $G$ on a \ca\  $A$.
Regard $K_* (A)$ as a $\Z [G]$-module via the action of $G$ on $A$.
Regard the equivariant K-theory $K_*^{\widehat{G}} (C^* (G, A, \af))$
of the dual action
${\widehat{\af}} \colon {\widehat{G}} \to \Aut (C^* (G, A, \af))$
as a $\Z [G]$-module as discussed above.
Then $K_*^{\widehat{G}} (C^* (G, A, \af)) \cong K_* (A)$
as $\Z [G]$-modules.
\end{lem}

\begin{proof}
Let $B = C^* ({\widehat{G}}, \, C^* (G, A, \af), \, {\widehat{\af}})$
be the second crossed product, with second dual action
$\bt \colon G \to \Aut (B)$.
Regard $K_* (B)$ as a $\Z [G]$-module via this action of $G$.
Proposition~2.7.10 of~\cite{PhT},
applied to the action ${\widehat{\af}}$ of ${\widehat{G}}$ on
$C^* (G, A, \af)$,
shows that $K_0^{\widehat{G}} (C^* (G, A, \af)) \cong K_0 (B)$
as $\Z [G]$-modules.
Applying this result to the suspension $S A$,
we obtain the same isomorphism for $K_1$ as well.
By Takai duality~\cite{Tk},
there is an isomorphism $B \cong A \otimes K (l^2 (G))$
which intertwines the action $\bt$
with the tensor product of $\af$ and an inner action on $K (l^2 (G))$.
It is therefore immediate that
$K_* (B) \cong K_* (A)$ as $\Z [G]$-modules.
\end{proof}

We immediately get a realization theorem using type~I \ca s.

\begin{thm}\label{RealizeEqKThT1}
Let $G_*$ be a countable $\Zt$-graded $R (\Zt)$-module.
Then there exists a separable type~I \ca\  $A$,
and an action $\af$ of $\Zt$ on $A$, such that
$K_*^{\Zt} (A) \cong G_*$ as $R (\Zt)$-modules.
\end{thm}

\begin{proof}
Identify $R (\Zt)$ with $\Z [\Zth]$
as discussed before Lemma~\ref{CrProd}.
Since $\Zth \cong \Zt$, we can
use Proposition~\ref{ExistAlg} to find a separable type~I \ca\  $B$
and an action $\bt \colon \Zth \to \Aut (B)$ such that
$K_* (B) \cong G_*$ as $R (\Zt)$-modules.
Let $A = \Cstd{B}{\bt}$, equipped with the dual
action $\af = {\widehat{\bt}}$ of $\Zt$.
Lemma~\ref{CrProd} implies that
$K_*^{\Zt} (A) \cong G_*$ as $R (\Zt)$-modules.
Clearly $A$ is separable,
and Theorem~4.1 of~\cite{RfFG} implies that $A$ is type~I\@.
\end{proof}

For the realization theorem using Kirchberg algebras,
we need another lemma.

\begin{lem}\label{DirLimIsK}
Let
\[
A_0 \stackrel{\rh_1}{\longrightarrow}
A_1 \stackrel{\rh_2}{\longrightarrow}
A_2 \stackrel{\rh_3}{\longrightarrow} \cdots
\]
be a direct system of \ca s in which each
$A_n$ is a finite direct sum of Kirchberg algebras
satisfying the \uct.
Suppose that for each $n$ the partial map determined by $\rh_n$ from
any summand of $A_{n - 1}$ to any other summand of $A_n$ is nonzero.
Then the direct limit of this system is again a Kirchberg algebra
satisfying the \uct.
\end{lem}

\begin{proof}
Let $A = \dirlim A_n$.
Clearly $A$ is separable and nuclear.
The condition that all the partial maps be nonzero, together with
the simplicity of all the summands, guarantees
that the algebraic direct limit of the $A_n$ is simple.
A standard argument now shows that $A$ is simple.
It is easy to check, using standard direct limit methods, that $A$ has
real rank zero and that every nonzero \pj\  in $A$ is properly infinite.
(For the second statement, use the fact that every \pj\  in $A$ is
equivalent to a \pj\  in the algebraic direct limit.)
Therefore $A$ is purely infinite simple.
Finally, $A$ satisfies the Universal Coefficient Theorem because
the class of such algebras is closed under direct limits.
This is essentially Proposition~2.3(b) of \cite{RS}, in view of
Theorem~4.1 of \cite{RS} and its converse (the converse being trivial).
\end{proof}

\begin{thm}\label{RealizeEqKTh}
Let $G_*$ be a countable $\Zt$-graded $R (\Zt)$-module.
Then there exists a unital Kirchberg algebra $A$,
and an action $\af$ of $\Zt$ on $A$, such that
$K_*^{\Zt} (A) \cong G_*$ as $R (\Zt)$-modules,
and such that $\Cst{A}{\af}$ is a Kirchberg algebra.
\end{thm}

\begin{proof}
Identify $R (\Zt)$ with $\Z [\Zth]$
as discussed before Lemma~\ref{CrProd}.
Let $\ta \in \Zth$ be the nontrivial element.

First suppose that
the action of $\ta$ on the $R (\Zt)$-module $G_*$ is nontrivial.
Since $\Zth \cong \Zt$, we can
use Theorem~\ref{LiftGpAction} to find a unital Kirchberg algebra $B$
and an action $\bt \colon \Zth \to \Aut (B)$ such that
$K_* (B) \cong G_*$ as $R (\Zt)$-modules.
Let $A = \Cstd{B}{\bt}$, equipped with the dual
action $\af = {\widehat{\bt}}$ of $\Zt$.
Lemma~\ref{CrProd} implies that
$K_*^{\Zt} (A) \cong G_*$ as $R (\Zt)$-modules.
Moreover, $\bt_{\ta}$ is outer (because it is nontrivial on K-theory),
so Corollary~4.4 of~\cite{JO}
implies that $A$ is a Kirchberg algebra.
Moreover $\Cst{A}{\af} \cong M_2 \otimes B$ is also a Kirchberg algebra.

Now suppose that
the action of $\ta$ on the $R (\Zt)$-module $G_*$ is trivial.
Choose a unital Kirchberg algebra $A_0$ satisfying the \uct\  %
such that $K_* (A_0) \cong G_*$ as abelian groups.
Let $\af^{(0)} \colon \Zt \to \Aut (A_0 \oplus A_0)$
be the action such that $\af^{(0)}_{\ta} (a, b) = (b, a)$.
We identify the crossed product
$C^* \left( \Zt, \, A_0 \oplus A_0, \, \af^{(0)} \right)$.
Let $u \in C^* \left( \Zt, \, A_0 \oplus A_0, \, \af^{(0)} \right)$
be the canonical unitary \oot.
Then there is an isomorphism
$\sm \colon C^* \left( \Zt, \, A_0 \oplus A_0, \, \af^{(0)} \right)
        \to M_2 (A_0)$,
determined by
\[
\sm (a, b)
 = \left( \begin{array}{cc} a & 0 \\ 0 & b \end{array} \right)
\]
for $a, \, b \in A_0$, and
\[
\sm (u)
 = \left( \begin{array}{cc} 0 & 1 \\ 1 & 0 \end{array} \right).
\]
The dual action is inner,
so $K_*^{\Zt} (A_0 \oplus A_0) \cong G_*$ as $R (\Zt)$-modules.

The algebra $A_0 \oplus A_0$ is not a Kirchberg algebra.
To remedy this, choose a nonzero \pj\  $q \in A_0$
such that $[q] = 0$ in $K_0 (A_0)$.
Set $p = 1 - q$, so that $[p] = [1_{A_0}]$ in $K_0 (A_0)$.
Choose
(\cite{Kr}; or Theorem~4.1.3 and the proofs
of Corollary~4.4.2 and Theorem~4.2.4 of~\cite{Ph})
unital \hm s $\ph \colon A_0 \to p A_0 p$ and
$\ps \colon A_0 \to q A_0 q$
such that
\[
A_0 \stackrel{\ph}{\longrightarrow} p A_0 p \longrightarrow A_0
\]
is the identity on K-theory and $\ps_* = 0$.
Define a unital equivariant \hm\  %
$\rh \colon A_0 \oplus A_0 \to A_0 \oplus A_0$
by $\ph (a, b) = (\ph (a) + \ps (b), \, \ph (b) + \ps (a))$
for $a, \, b \in A_0$.
Let $A$ be the direct limit of the system
\[
A_0 \oplus A_0 \stackrel{\rh}{\longrightarrow}
 A_0 \oplus A_0 \stackrel{\rh}{\longrightarrow}
 A_0 \oplus A_0 \stackrel{\rh}{\longrightarrow} \cdots.
\]
Then $A$ is a unital Kirchberg algebra by Lemma~\ref{DirLimIsK}.
Let $\af \colon \Zt \to \Aut (A)$ be the direct limit action.

We next show that $K_*^{\Zt} (A) \cong G_*$ as $R (\Zt)$-modules.
We do this by proving that
$\rh_* \colon K_*^{\Zt} (A_0 \oplus A_0) \to K_*^{\Zt} (A_0 \oplus A_0)$
is an isomorphism
and using Proposition~2.5.4 of~\cite{PhT}
on the direct system and its suspension.
We claim that the corresponding \hm\  %
\[
{\overline{\rh}} \colon
 {\textstyle{ C^* \left( \Zt, \, A_0 \oplus A_0, \, \af^{(0)} \right) }}
\to
 {\textstyle{ C^* \left( \Zt, \, A_0 \oplus A_0, \, \af^{(0)} \right) }}
\]
is, under the identification of the crossed product above,
given by
\[
{\overline{\rh}} (x)
 = (\id_{M_2} \otimes \ph) (x) + u (\id_{M_2} \otimes \ps) (x) u^*.
\]
It is enough to check this on the image of
$A_0 \oplus A_0$, namely the diagonal matrices, and on $u$;
this is easy.
It is immediate from the choice of $\ph$ and $\ps$ that this
map is the identity on K-theory.

It remains to show that $\Cst{A}{\af}$ is a Kirchberg algebra.
This algebra is the direct limit of the system of crossed products
\[
C^* ( \Zt, \, A_0 \oplus A_0, \, \af^{(0)} )
     \stackrel{ {\overline{\rh}} }{\longrightarrow}
 C^* ( \Zt, \, A_0 \oplus A_0, \, \af^{(0)} )
     \stackrel{ {\overline{\rh}} }{\longrightarrow}
 C^* ( \Zt, \, A_0 \oplus A_0, \, \af^{(0)} )
     \stackrel{ {\overline{\rh}} }{\longrightarrow} \cdots,
\]
which, by the computation at the beginning of this case,
can be rewritten as
\[
M_2 (A_0) \longrightarrow
 M_2 (A_0) \longrightarrow
 M_2 (A_0) \longrightarrow \cdots.
\]
The direct limit of this system
is a Kirchberg algebra by Lemma~\ref{DirLimIsK}.
\end{proof}

\begin{rmk}\label{NucQ}
When the action of the nontrivial element of $\Zth$ is nontrivial,
we do not know whether
the \ca\  in Theorem~\ref{RealizeEqKTh} satisfies the \uct.
\end{rmk}

\section{An explicit example}\label{Sec:ExplEx}

In this section, we write down a formula for an automorphism
in the smallest nontrivial case of Theorem~\ref{LiftGpAction}.
This is the case $K_0 (A) \cong \Zq{3}$, $K_1 (A) = 0$,
and $\gm_0$ is the unique nontrivial automorphism of $\Zq{3}$.
Since $\gm_0$ fixes only the identity element of $\Zq{3}$,
the automorphism will exist only when $[1_A] = 0$ in $K_0 (A)$.
Let ${\mathcal{O}}_{4}$ be the Cuntz algebra~\cite{Cu1},
with K-theory as computed in~\cite{Cu2}.
Then the unital Kirchberg algebra $A$ with this K-theory
and satisfying the \uct\  is $M_3 \otimes {\mathcal{O}}_{4}$.

\begin{exa}\label{M4O3}
We give an explicit formula, in terms of the standard generators,
for an \ota\  $\ph$ of
$A = M_3 \otimes {\mathcal{O}}_{4}$
which is nontrivial on $K_0 (A)$.

Let $s_1$, $s_2$, $s_3$, and $s_4$ be the standard generating
isometries of ${\mathcal{O}}_{4}$,
and define $p_m = s_m s_m^*$.
Further let $e_{j, k}$, for $1 \leq j, \, k \leq 3$,
be the standard matrix units of $M_3$,
satisfying $e_{j, k} e_{k, l} = e_{j, l}$ etc.
Then $\ph$ is determined by the formulas:
\begin{align*}
e_{1, 1} \otimes 1 & \mapsto
      f_{1, 1} = (e_{2, 2} + e_{3, 3}) \otimes 1  \\
e_{2, 2} \otimes 1 & \mapsto
      f_{2, 2} = e_{1, 1} \otimes (p_1 + p_2)  \\
e_{3, 3} \otimes 1 & \mapsto
      f_{3, 3} = e_{1, 1} \otimes (p_3 + p_4)  \\
e_{1, 2} \otimes 1 & \mapsto
      f_{1, 2} = e_{2, 1} \otimes s_1^* + e_{3, 1} \otimes s_2^*  \\
e_{1, 3} \otimes 1 & \mapsto
      f_{1, 3} = e_{2, 1} \otimes s_3^* + e_{3, 1} \otimes s_4^*  \\
e_{1, 1} \otimes s_1 & \mapsto
      v_1 = e_{2, 2} \otimes s_1 + e_{2, 3} \otimes s_2  \\
e_{1, 1} \otimes s_2 & \mapsto
      v_2 = e_{2, 2} \otimes s_3 + e_{2, 3} \otimes s_4  \\
e_{1, 1} \otimes s_3 & \mapsto
      v_3 = e_{3, 2} \otimes s_1 + e_{3, 3} \otimes s_2  \\
e_{1, 1} \otimes s_4 & \mapsto
      v_4 = e_{3, 2} \otimes s_3 + e_{3, 3} \otimes s_4
\end{align*}

One must check two things:
that $\ph$ is a \hm, and that $\ph^2 = \id_A$.
For the first, one checks the following relations:
\begin{itemize}
\item
$f_{1, 1}$, $f_{2, 2}$, and $f_{3, 3}$
are orthogonal \pj s which sum to $1$.
\item
$f_{1, j} f_{1, j}^* = f_{1, 1}$ and $f_{1, j}^* f_{1, j} = f_{j, j}$
for $j = 2, \, 3$.
\item
$v_m^* v_m = f_{1, 1}$ for $1 \leq m \leq 4$ and
$\sum_{m = 1}^4 v_m v_m^* = f_{1, 1}$.
\end{itemize}
The details of the computation are somewhat long, and are omitted.

To prove that $\ph^2 = \id_A$, one checks that $\ph (\ph (a)) = a$
for each of the generators used above.
Again, we omit the details.
\end{exa}

\end{document}